\documentclass[12pt]{article}

\usepackage{amsmath}
\usepackage{amssymb}
\usepackage{amsthm}
\usepackage{fullpage}
\usepackage{hyperref}

\title{A note on super Catalan numbers}
\author{ Evangelos Georgiadis\thanks{Massachusetts Institute of Technology, Cambridge, M.A. 02139, U.S.A.  \textit{Email:}~\texttt{egeorg@mit.edu}}\and Akihiro Munemasa\thanks{Graduate School of Information Sciences, Tohoku University, Sendai 980-8579, Japan. \textit{Email:}~\texttt{munemasa@math.is.tohoku.ac.jp}}\and Hajime Tanaka\thanks{Graduate School of Information Sciences, Tohoku University, Sendai 980-8579, Japan. \textit{Email:}~\texttt{htanaka@m.tohoku.ac.jp}}}
\date{May 22, 2012}

\begin{document}

\maketitle

\begin{abstract}
We show that the super Catalan numbers are special values of
the Krawtchouk polynomials by deriving an expression for
the super Catalan numbers in terms of a signed set.
\end{abstract}

\bigskip

The super Catalan numbers,
\begin{equation*}
	S(m,n) := \frac{(2m)!(2n)!}{m!n!(m+n)!},
\end{equation*}
as designated by Gessel in ~\cite[Eq.~(28)]{Gessel1992JSC} form, hierarchically speaking, special cases of super ballot numbers (cf.~\cite[pp.~180, 189]{Gessel1992JSC}).
Historically, Gessel points out that these numbers had been observed as early as 1874 and studied by E. Catalan~\cite{Catalan1874JEPN}; Aguiar and Hsiao in~\cite{Aguiaretal2005EJC} provide a more detailed account of earlier appearences (cf.~\cite{Catalan1873}, \cite{Catalan1883},~\cite{Catalan1887} and Riordan~\cite[Chapter 3, Exercise 9, p.~120]{Riordan1968}).
Surprisingly, the innocent looking numbers, $S(m,n)$, seem to have been immune against a combinatorial interpretation for all values of $(m,n)$ over the last century despite limited success stories for particular values (see problem 66(a) in Stanley's bijective open problems compendium~\cite{StanleyOP}).
The following references highlight the success cases.
In~\cite{Gessel1992JSC}, Gessel notes that for $S(1,n)/2$ we obtain the Catalan number $C_n$; whereas for the case when $m=0$, we yield middle binomial coefficients, $\binom{2n}{n}$.
In~\cite{Gesseletal2005JIS}, Gessel and Xin provide a combinatorial interpretation in terms of Dyck paths when $m=2$ or 3.
An alternative combinatorial interpretation for the case $m=2$ was provided by Schaeffer in~\cite{Schaeffer2003} using a method that was introduced in the interpretation to formulas of Tutte for planar maps.
A more topologically flavored yet still combinatorial interpretation for the $m=2$ case is also available by Pippenger and Schleich in~\cite{Pippengeretal2003}; they count cubic trees on $n$ interior vertices (or the number of hexagonal trees with $n$ nodes).
In 2005, Callan in~\cite{Callan2005JIS} provided an elegant combinatorial interpretation of the recurrence $S(m,n)/2=\sum_{k \geq 0} 2^{n-m-2k} \binom{n-m}{2k}S(m,k)/2$ for the case when $m=2$ showing that it enumerates the aligned cubic trees by number of vertices that are neither a leaf nor adjacent to a leaf.

In this note we establish the following expression for super Catalan numbers:
\begin{equation}\label{S1}
	S(m,n)=(-1)^m\sum_{P\in\mathcal{P}_{m+n}}(-1)^{h_{2m}(P)},
\end{equation}
where the sum is over the set $\mathcal{P}_{m+n}$ of all lattice paths from $(0,0)$ to $(m+n,m+n)$ consisting of unit steps to the right and up, and $h_{2m}(P)$ denotes the height of $P=(P_0,P_1,\dots,P_{2(m+n)})\in\mathcal{P}_{m+n}$ after the $2m^{\mathrm{th}}$ step, i.e., the $y$-coordinate of $P_{2m}$.
Although
this is an interpretation of $S(m,n)$ in terms of a signed set only,
the right-hand side of \eqref{S1} is a special value of the
Krawtchouk polynomial defined as follows:
\begin{equation*}
	K_j^d(x)=\sum_{h=0}^j(-1)^h\binom{x}{h}\binom{d-x}{j-h}.
\end{equation*}
Then \eqref{S1} is equivalent to
\begin{equation}\label{S2}
	K_{m+n}^{2(m+n)}(2m)=(-1)^m S(m,n).
\end{equation}
To see the equivalence, 
observe that each $P\in\mathcal{P}_{m+n}$ has exactly $m+n$ up-steps and that the number of $P\in\mathcal{P}_{m+n}$ with $h_{2m}(P)=h$ is therefore equal to $\binom{2m}{h}\binom{2n}{m+n-h}$.

Krawtchouk polynomials $K_j^d(x)$ appear as the coefficients of the so-called
MacWilliams identities (cf.~\cite[Chap.~5, \S 2]{MS}), and
also as the eigenvalues of the distance-$j$ graph of the $d$-cube
(cf.~\cite[Chap.~3, \S 2]{BannaiIto1984}).
The identity shows that $\{(-1)^mS(m,n)\mid m,n\geq0,\;m+n=N\}$ coincides with the set of non-zero eigenvalues of the distance-$N$ graph of the $2N$-cube, which is known as
the orthogonality graph and has been studied in connection with pseudo-telepathy in quantum information theory (cf.~\cite{GN2008SIAM}).
Finally, \eqref{S2} follows immediately from the identity of von Szily (cf.~\cite[Eq.~(29)]{Gessel1992JSC}):
\begin{align*}
	S(m,n)&=\sum_{k\in\mathbb{Z}}(-1)^k\binom{2m}{m+k}\binom{2n}{n-k} \\
	&=(-1)^m\sum_{h=0}^{m+n}(-1)^h\binom{2m}{h}\binom{2n}{m+n-h} \\
	&=(-1)^mK_{m+n}^{2(m+n)}(2m).
\end{align*}
We note that \eqref{S2}, as well as the identity of von Szily, is just a restatement of (a special case of) Kummer's evaluation of well-poised $_2F_1(-1)$ series.

To obtain a proper interpretation as the size of a set of certain paths, 
we need to find an injection from the set
\[
\{P\in\mathcal{P}_{m+n}\mid h_{2m}(P)\not\equiv m\pmod{2}\}
\]
to
\[
\{P\in\mathcal{P}_{m+n}\mid h_{2m}(P)\equiv m\pmod{2}\},
\]
and a description of the complement of the image. This is known for
the case $m=1$ (see \cite[Section~5.3]{Ai}), 
but it seems to be a difficult problem in general.

\section*{Acknowledgments}
Special thanks are due to Ole Warnaar and Richard Askey for valuable comments. E.G.~is indebted to M. Sipser and J. Tsitsiklis of MIT and would like to thank I. Gessel of Brandeis University for helpful background comments on this problem in early 2010.
H.T.~was supported in part by the JSPS Excellent Young Researchers Overseas Visit Program.

\end{document}